\documentclass{article}

\usepackage{arxiv}

\usepackage[utf8]{inputenc} % allow utf-8 input
\usepackage[T1]{fontenc}    % use 8-bit T1 fonts
\usepackage{hyperref}       % hyperlinks
\usepackage{url}            % simple URL typesetting
\usepackage{booktabs}       % professional-quality tables
\usepackage{amsfonts}       % blackboard math symbols
\usepackage{nicefrac}       % compact symbols for 1/2, etc.
\usepackage{microtype}      % microtypography
\usepackage{lipsum}

\usepackage{graphicx}
\usepackage{amsmath}

\title{Fast Computation of the Direct Scattering Transform by Fourth Order Conservative Multi-Exponential Scheme}

\author{
  Sergey Medvedev$^{1,2,*}$, Igor Chekhovskoy$^{2,1}$, Irina Vaseva$^{1,2}$, Mikhail Fedoruk$^{2,1}$\\
$^{1}$ Institute of Computational Technologies, SB RAS, Novosibirsk
630090, Russia,\\
$^{2}$ Novosibirsk State University, Novosibirsk 630090, Russia,\\
* Corresponding author: medvedev@ict.nsc.ru
}

\begin{document}
\maketitle

\begin{abstract}
A fourth-order multi-exponential scheme is proposed for the Zakharov-Shabat system. The scheme represents a product of 13 exponential operators. The construction of the scheme is based on a fourth-order three-exponential scheme, which contains only one exponent with a spectral parameter. This exponent is factorized to the fourth-order with the Suzuki formula of 11 exponents. The obtained scheme allows the use of a fast algorithm in calculating the initial problem for a large number of spectral parameters and conserves the quadratic invariant exactly for real spectral parameters.
\end{abstract}

% keywords can be removed

\section{Introduction}
In 1971, Zakharov and Shabat (ZS) showed that the nonlinear Schr\"odinger equation (NLSE)
\begin{equation}\label{nlse}
i\frac{\partial q}{\partial z}+\frac{\sigma}{2}\frac{\partial^2 q}{\partial t^2}+|q|^2q=0,\quad \sigma=\pm 1
\end{equation}
can be integrated by the inverse problem method (or so-called nonlinear Fourier transform -- NFT) previously applied to the Korteweg de Vries equation~\cite{ZakharovShabat1972}. After that, interest in the NLSE arose in all areas of physics connected with wave systems, because the NLSE describes the envelope for narrow wave beams.
In 1973, Hasegawa and Tappert numerically investigated the NLSE in respect to the propagation of light pulses in optical fibers~\cite{Hasegawa1973a}. They proposed using solitons as an information carrier for fiber lines with anomalous dispersion at $\sigma=1$. For normal dispersion at $\sigma =-1$, solitons do not exist, as is well known.

Since that time, the study of the NLSE and its generalizations to describe the propagation of light pulses in optical fibers has begun. Analytical and numerical studies were carried out, as well as work on the development of numerical methods for integrating NLSE. Currently, the most popular and effective method is splitting into physical processes, the so-called split-step Fourier method (SSFM)~\cite{Hardin_SIAM_1973}.

On the other hand, attempts to create fast numerical algorithms for solving the inverse scattering problem for the NLSE have not stopped. Such methods are combined under the general name Fast Nonlinear Fourier Transform (FNFT)~\cite{Wahls2013, Wahls2015a, Wahls2016, Turitsyn2017Optica}.

In this paper, we propose a special fourth-order numerical method for solving the direct ZS problem (ZSP) and a fast algorithm for its numerical implementation.
The main advantage of the presented scheme is the conservation of the quadratic invariant for real spectral parameters, even in the fast version. This is the first proposed fast scheme with such property for the best of our knowledge.
The quadratic invariant conservation by numerical scheme allows calculating precisely the reflection coefficient, what is valuable for various telecommunication problems connected with NFT-based coding schemes (for example, NFDM~\cite{Yousefi2014III} and $b$-modulation~\cite{Wahls2017}).

\section{Direct spectral ZS problem}\label{sec:headings}

Direct spectral ZS problem for the NLSE~(\ref{nlse}) with the complex spectral parameter $\zeta$ can be rewritten as an evolutionary system
\begin{equation}\label{psit}
\frac{d{\Psi}(t)}{dt}=Q(t){\Psi}(t),
\end{equation}
where $q=q(t,z_0)$ is the initial field for the NLSE at the point $z_0$, which is the potential in the ZS problem, and
$$
{\Psi}(t)=\left(\begin{array}{c}\psi_1(t)\\\psi_2(t)\end{array}\right),\quad
Q(t)=\left(\begin{array}{cc}-i\zeta&q\\-\sigma q^*&i\zeta\end{array}\right).
$$
Here $z_0$ plays the role of the parameter and we will omit it. For details, we refer to the numerous literature, in particular,~\cite{medvedev2019exponential}.

Moreover, the system~(\ref{psit}) can be written in the gradient form
\begin{equation}
\left(\begin{array}{c}\psi_1\\\psi_2\end{array}\right)_t=KD
\left(\begin{array}{c}\psi_1\\\psi_2\end{array}\right)=K
\left(\begin{array}{c}\frac{\partial H}{\partial\psi_1^*}\\\frac{\partial H}{\partial\psi_2^*}\end{array}\right),
\end{equation}
where $H=({\Psi}^*,D{\Psi})$,
\begin{equation}
K=\left(\begin{array}{cc}-i\zeta&\sigma q\\-\sigma q^*&i\sigma\zeta\end{array}\right),\quad 
D=\left(\begin{array}{cc}1&0\\0&\sigma\end{array}\right).
\end{equation}
For real~$\zeta=\xi$ the matrix~$K$ before the gradient becomes anti-Hermitian $K=-K^\dagger$ for any $\sigma=\pm 1$ and, therefore, the system~(\ref{psit}) will conserve the quantity~$H$.

Assuming that $q(t)$ decays rapidly when $t\to \pm \infty $, the specific solutions (Jost functions) for ZSP~(\ref{psit}) can be derived as
\begin{equation}\label{psi0}
\Psi =\left(
\begin{array}{c}
\psi_{1}\\\psi_{2}
\end{array}
\right) = \left(
\begin{array}{c}
e^{-i\zeta t}\\0
\end{array}
\right)[1+o(1)],\quad t\to-\infty,
\end{equation}
and
\begin{equation}\label{psi0right}
\Phi =\left(
\begin{array}{c}
\phi_{1}\\\phi_{2}
\end{array}
\right) = \left(
\begin{array}{c}
0\\e^{i\zeta t}
\end{array}
\right)[1+o(1)],\quad t\to \infty,
\end{equation}
Then we obtain the Jost scattering coefficients $a(\xi )$ and $b(\xi )$ as
\begin{equation}\label{ab}
a(\xi)=\lim_{t\to\infty}\,\psi_1(t,\xi)\,e^{i\xi t},\quad b(\xi)=\lim_{t\to\infty}\,\psi_2(t,\xi)\,e^{-i\xi t}.
\end{equation}
The functions $a(\xi)$ and $b(\xi)$ can be extended to the upper half-plane $\xi \to \zeta $, where $\zeta $ is a complex number with the positive imaginary part $\eta =Im\zeta >0$. The spectral data of ZSP~(\ref{psit}) are determined by $a(\zeta )$ and $b(\zeta )$ in the following way:\\
(1) zeros of $a(\zeta)=0$ define the discrete spectrum $\{\zeta_k\}$, $k=1,...,K$ of ZSP~(\ref{psit}) and phase coefficients
$$r_k=\left.\frac{b(\zeta)}{a'(\zeta)}\right|_{\zeta=\zeta_k},\quad\mbox{where}\quad a'(\zeta)=\frac{da(\zeta)}{d\zeta};$$
(2) the continuous spectrum is determined by the reflection coefficient
$$r(\xi)={b(\xi)}/{a(\xi)},\quad \xi\in\mathbb{R}. $$

These spectral data were defined using the "left"\,  boundary condition (\ref{psi0}). Both conditions (\ref{psi0}) and (\ref{psi0right}) can be used to calculate the coefficient~$b(\zeta _{k} )$ of the discrete spectrum with the relation $\Psi(t,\zeta_k)=\Phi(t,\zeta_k)b(\zeta_k)$.

For real values of the spectral parameter $\zeta=\xi$, we have the quadratic invariant~$H=|\psi_1|^2+\sigma|\psi_2|^2$.
Taking into account the boundary conditions~(\ref{psi0}), we get the same condition $H=1$ for $\sigma=\pm 1$.

\section{Computational features of ZS system}

We solve a linear system of the form~(\ref{psit}) with the matrix~$Q(t)$ linearly dependent on the complex function~$q(t)$. The numerical implementation of the continuous function~$ q(t)$ is a discrete function~$q_n=q(t_n)$, which is defined at the integer nodes $t_n$ of the uniform grid with the step $\tau$. Since we are considering a finite time interval, we will solve the problem on the interval $[-L,L]$ with the total number of points equal to $M+1$, the grid step in this case is $\tau=2L/M$ and $t_n=-L+\tau n$, where $n=0,...,M$.

The main features of the discrete problem are the following:
\begin{enumerate}
\item The matrix of the system $Q$ is given on a uniform grid with a step $\tau$, therefore the unknown function $\Psi$ must also be calculated on a uniform grid with a step $\tau$. One cannot use Runge-Kutta methods on such grid. If, for example, an explicit 4th-order RC scheme is used, then it is necessary to take the computational grid with a double step. In this case, the values of $Q_n$ will enter unequally.
\item For small values of the potential $|q(t)|<<|\zeta|$ and $\mbox{Im}\,\zeta>0$, the ZS system contains exponentially increasing and decreasing solutions therefore A-stability of difference methods is required~\cite{dahlquist1963special}. {\it Dalquist's Second Barrier} restricts multistep methods.This barrier says that there are no explicit A-stable multistep methods. 2nd order of convergence is maximum for implicit multi-step methods~\cite{hairer1991solving}.
\item The ZS system has a second-order matrix; therefore, inverse matrices and the exponent of the matrix are easily calculated. This allows us to include almost any functions of the matrix~$Q$ in difference schemes.
\item To calculate the spectral data, it is required to solve the ZS system for a large number of spectral parameter values~$\zeta$ with the fixed potential $q(t)$. If possible, this should be taken into account when implementing the algorithms.
\end{enumerate}

\section{Scheme}

Previously, we found the necessary conditions for the existence of a one-step fourth-order scheme $\Psi(t_n+\tau/2)=T_n\Psi(t_n-\tau/2)$ in the form of the Taylor series for the transition matrix~$T_n$. We obtained several fourth-order schemes that exactly conserve the quadratic invariant~$H$. However, only one scheme in the form of three exponentials is convenient for the fast computation~\cite{medvedev2019exponential}:
\begin{equation}\label{tripleexp}
T_n=e^{\left\{\frac{\tau^2}{12}Q^{(1)}_n+\frac{\tau^3}{48}Q^{(2)}_n\right\}}e^{\tau Q_n}e^{\left\{-\frac{\tau^2}{12}Q^{(1)}_n+\frac{\tau^3}{48}Q^{(2)}_n\right\}},
\end{equation}
here $Q_n=Q(t_n)$ and $Q^{(k)}_n$, $k =1,2,$ are the $k$-th derivative of the matrix $Q$ approximated by central finite differences of the second order.

To construct the fast algorithm, we must express $\exp(\tau Q)$, where $Q=A+B$, and
\begin{equation}
A=\left(\begin{array}{cc}-i\zeta&0\\0&i\zeta\end{array}\right),\quad
B(t)=\left(\begin{array}{cc}0&q\\-\sigma q^*&0\end{array}\right),
\end{equation}
in the form of a polynomial in exponentials of $A$ and $B$ with rational weights. For example, one can use the expansions for $\exp(\tau Q)$ suggested in \cite{Prins2018a}. However, representing the sum of the product of exponentials does not guarantee the exact conservation of the invariant~$H$. In order for the scheme to be suitable for the fast algorithm and conserve the invariant~$H$, it suffices to represent the matrix $\exp(\tau Q)$ as the product of exponentials of $A$ and $B$ with real rational coefficients.
Since for $\sigma=1$ the matrices $A$ and $B$ are Hermitian, then, in this case, each exponent will be unitary and the resulting scheme will conserve the quadratic invariant~$H$. Conserving~$H$ also takes place for $\sigma=-1$. Details can be found in \cite{medvedev2019exponential}. The rationality condition for weight coefficients provides an opportunity to represent the transition matrix in the form of the ratio of two polynomials in~$\exp(A)$.

\section{Suzuki factorization}
Since the scheme~(\ref{tripleexp}) has a fourth order of accuracy in $\tau$, it is necessary to have factorization of the same order. In addition, factorization should be suitable for a fast algorithm, i.e. have rational ratios. An example of such factorization is given in
\cite{suzuki1992general}:
\begin{equation}\label{factor4}
e^{\tau(A+B)}=
e^{\frac{7}{48}\tau B}e^{\frac{1}{3}\tau A}e^{\frac{3}{8}\tau B}e^{-\frac{1}{3}\tau A}e^{-\frac{1}{48}\tau B} e^{\tau A}e^{-\frac{1}{48}\tau B}e^{-\frac{1}{3}\tau A}
e^{\frac{3}{8}\tau B}e^{\frac{1}{3}\tau A}e^{\frac{7}{48}\tau B}.\nonumber
\end{equation}
We introduce the notation $Z = \exp(-\frac{i}{3}\tau\zeta)$, then the three exponents participating in this expansion take the form
\begin{equation}
e^{\frac{1}{3}\tau A}=\left(\begin{array}{cc}Z&0\\0&Z^{-1}\end{array}\right)=Z^{-1}\left(\begin{array}{cc}Z^2&0\\0&1\end{array}\right),
\end{equation}
\begin{equation}
e^{-\frac{1}{3}\tau A}=\left(\begin{array}{cc}Z^{-1}&0\\0&Z\end{array}\right)=Z^{-1}\left(\begin{array}{cc}1&0\\0&Z^2\end{array}\right),
\end{equation}
\begin{equation}
e^{\tau A}=\left(\begin{array}{cc}Z^3&0\\0&Z^{-3}\end{array}\right)=Z^{-3}\left(\begin{array}{cc}Z^6&0\\0&1\end{array}\right).
\end{equation}
Thus, the right-hand side of Eq.~(\ref{factor4}) is a rational function
\begin{equation}\label{SZ}
\frac{S(Z)}{Z^7},
\end{equation}
where $S(Z)$ is a polynomial not higher than 14 degrees in $Z$.
Since $Z$ is included only in the square, it is possible to introduce the variable $W=Z^2$, then (\ref{SZ}) takes the form
\begin{equation}\label{SW}
\frac{\hat{S}(W)}{W^{\frac{7}{2}}},
\end{equation}
where $\hat{S}(W)$ is a polynomial of degree 7 or less in $W$. The denominator is taken out and calculated independently.
It should be noted, that except of factorization~(\ref{factor4}) the symmetrical representation can be applied, when the matrices $A$ and $B$ are interchanged.
Such factorization leads to the polynomial~$\hat{S}(W)$ of degree 52, which is more computationally difficult and is less accurate, thus it is not considered here.

\section{Numerical examples}

The presented scheme implementation was based on FNFT software library~\cite{Wahls2018}. It was compared with
Boffetta-Osborne scheme (BO)~\cite{Boffetta1992a} and
the triple-exponential scheme with not factorized exponential~(\ref{tripleexp}) (TES4).
All these algorithms conserve the quadratic invariant~$H$ for real spectral parameters~$\xi$.
We do not compare our scheme with other fast algorithms (see, for example, \cite{Prins2018a, Wahls2018, Vaibhav2018}), due to the theoretical lack of such property among them.
We have considered both variants of the scheme with Suzuki factorization: conventional (TES4SB) and fast (FTES4SB). The last letter in the scheme name denotes the decomposition type: TES4SA denotes the scheme with the exponential with matrix~$A$ at the edges of Suzuki decomposition, while TES4SB is referred to the scheme~(\ref{factor4}).

A model signal was considered in the form of a chirped hyperbolic secant
\begin{equation}\label{Chirped}
q(t) = A[\mbox{sech}(t)]^{1+iC}.
\end{equation}
with the following parameters: $A = 5.2$, $C = 4$ for both anomalous and normal dispersion.
The detailed analytical expressions of the spectral data for this type of potentials can be found in~\cite{medvedev2019exponential}.

\begin{figure}[htbp]
\centering
\includegraphics[width=0.6\linewidth]{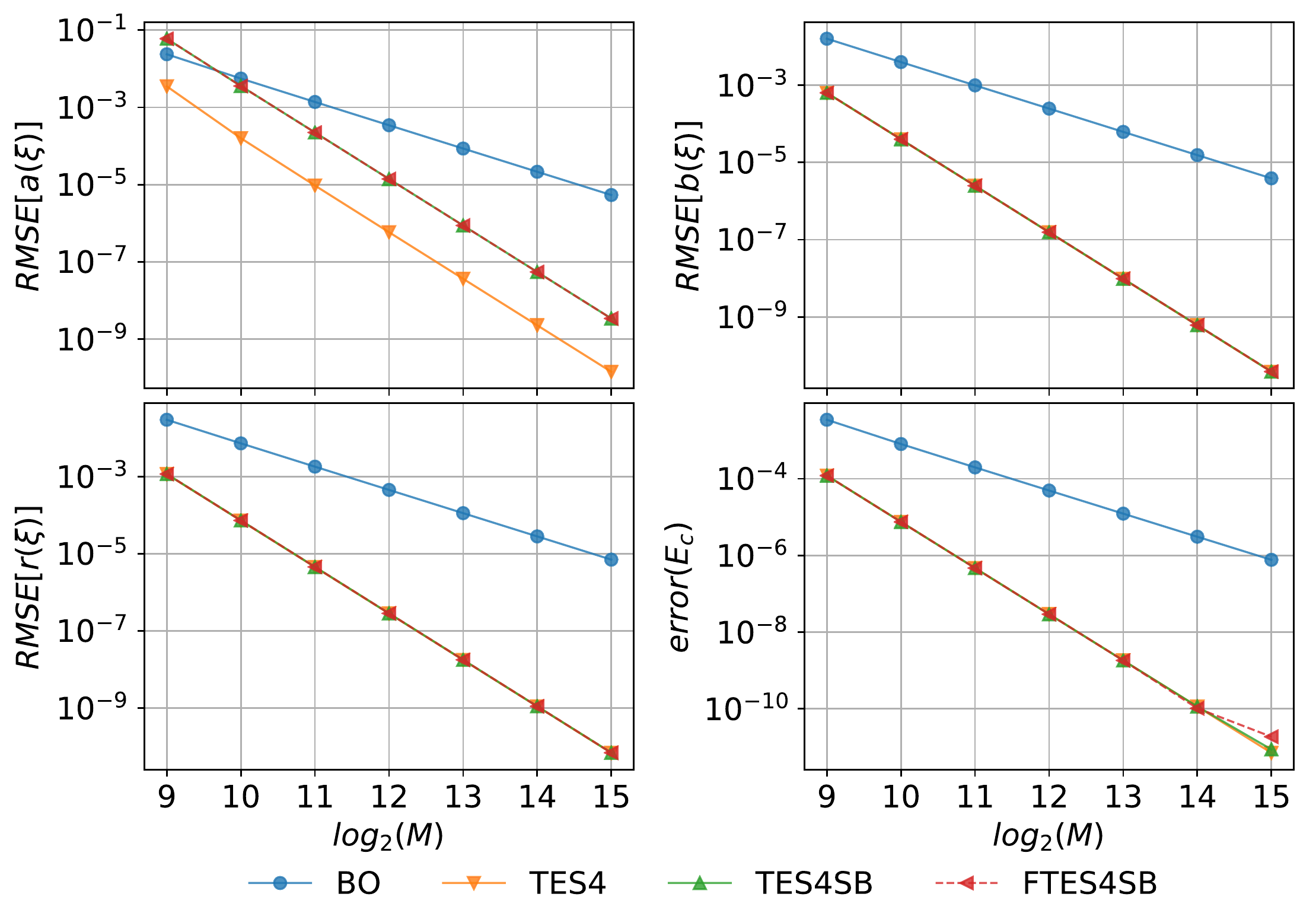}
\caption{The continuous spectrum errors for the chirped hyperbolic secant~(\ref{Chirped}) in the case of anomalous dispersion~$\sigma=1$.}
\label{fig:errorS1}
\end{figure}

\begin{figure}[htbp]
\centering
\includegraphics[width=0.6\linewidth]{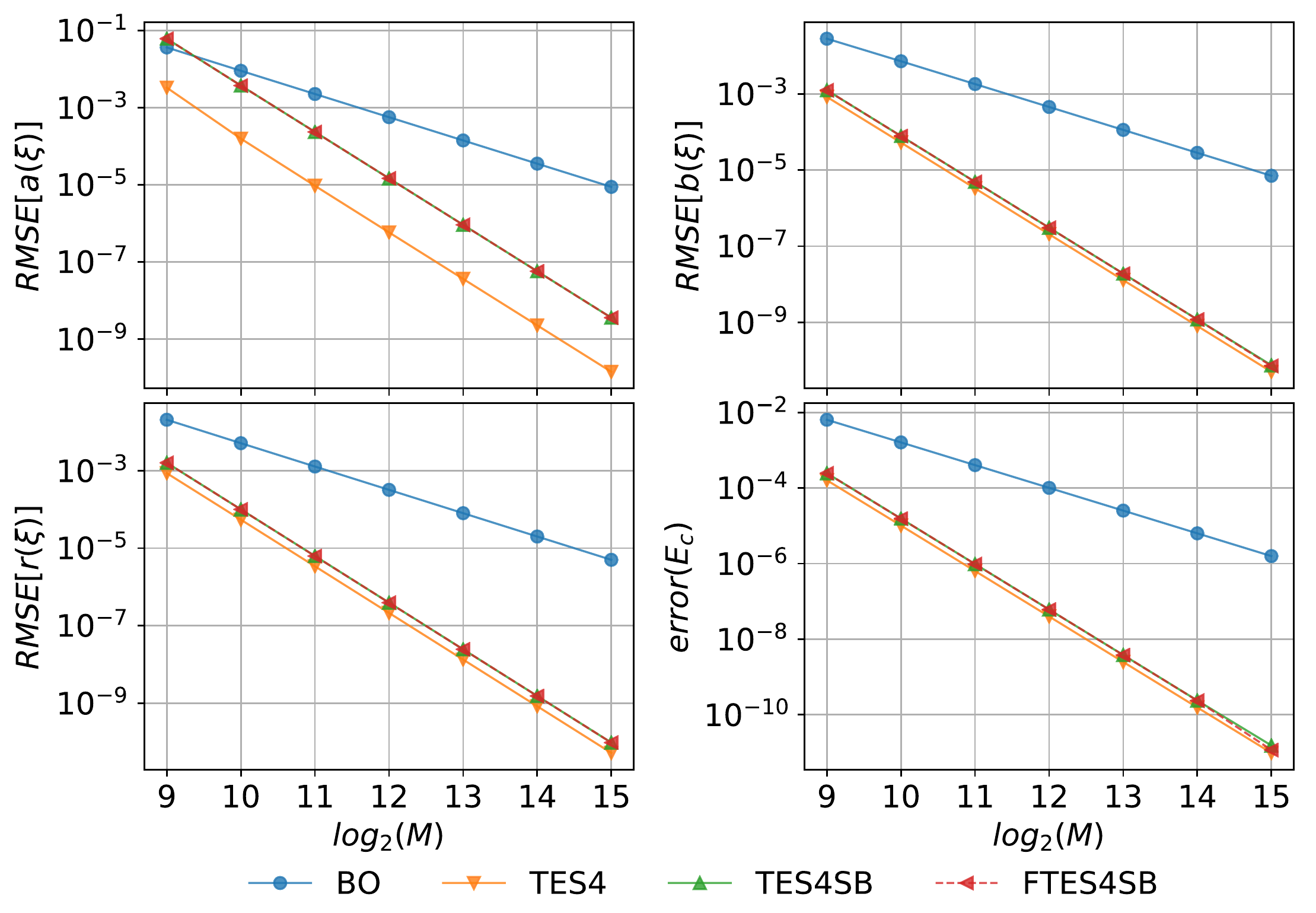}
\caption{The continuous spectrum errors for the chirped hyperbolic secant~(\ref{Chirped}) in the case of normal dispersion~$\sigma=-1$.}
\label{fig:errorSm1}
\end{figure}

We present the numerical errors of calculating the spectral data for continuous spectrum only, because of focusing on the conservation of the invariant~$H$ for real spectral parameters~$\xi$. To find the calculation errors of the continuous spectrum energy~$E_c$ the following formula was used
\begin{equation}\label{error}
\mbox{error}[E_c]=\frac{|E_c^{comp} - E_c^{exact}|}{\phi_0}, \quad
\phi_0 = 
\begin{cases}
E_c^{exact}, \mbox{ if } |E_c^{exact}| > 1\\
1, \mbox{ otherwise},
\end{cases}
\end{equation}
For the continuous spectrum we calculated the root mean squared error
\begin{equation}\label{MSE}
RMSE[\phi]=\sqrt{\frac{1}{N}\sum_{j=1}^{N}\frac{|\phi^{comp}(\xi_j) - \phi^{exact}(\xi_j)|^2}{|\phi_0(\xi_j)|^2}}, \quad
\phi_0 = 
\begin{cases}
\phi^{exact}(\xi_j), \mbox{ if } |\phi^{exact}(\xi_j)| > 1\\
1, \mbox{ otherwise},
\end{cases}
\end{equation}
where $\phi$ can represent $a(\xi)$, $b(\xi)$, $r(\xi)$ or $|H^{comp}(\xi) - H^{exact}(\xi)|$. Here we assume the spectral parameter $\xi \in [-20, 20]$ with the total number of points $N = 1025$.

Figures~\ref{fig:errorS1}, \ref{fig:errorSm1} present the errors calculated using the schemes under consideration. All schemes with triple-exponential representation demonstrated 4th order for the model signal and similar values of error in the case of anomalous dispersion.
%And only the scheme FNFT\_TES4SA showed a loss of the approximation order for normal dispersion.

\begin{figure}[htbp]
\centering
\includegraphics[width=0.6\linewidth]{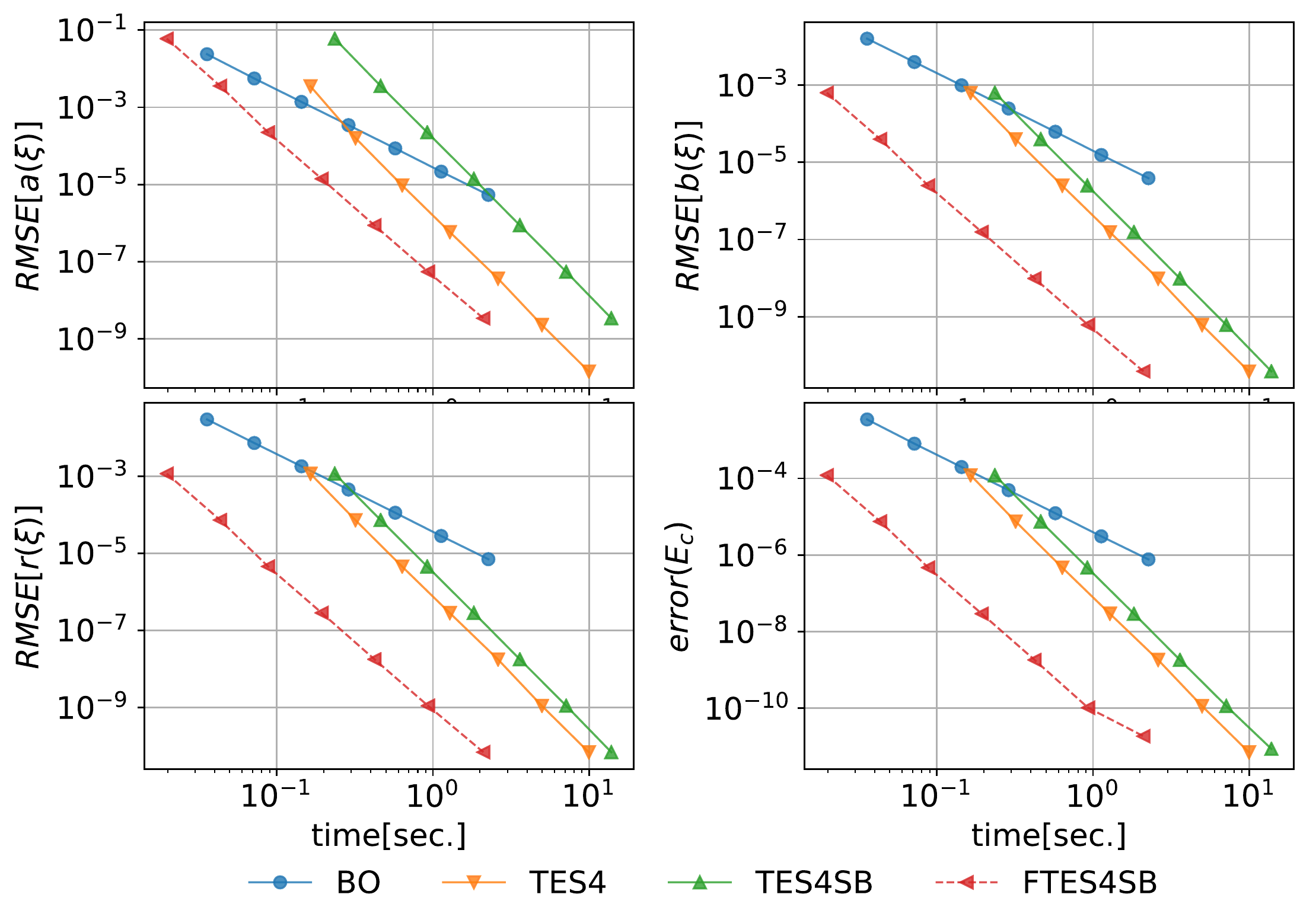}
\caption{The continuous spectrum errors depending on the execution time trade-off for the chirped hyperbolic secant~(\ref{Chirped}) in the case of anomalous dispersion~$\sigma=1$.}
\label{fig:errorTimeS1}
\end{figure}

\begin{figure}[htbp]
\centering
\includegraphics[width=0.6\linewidth]{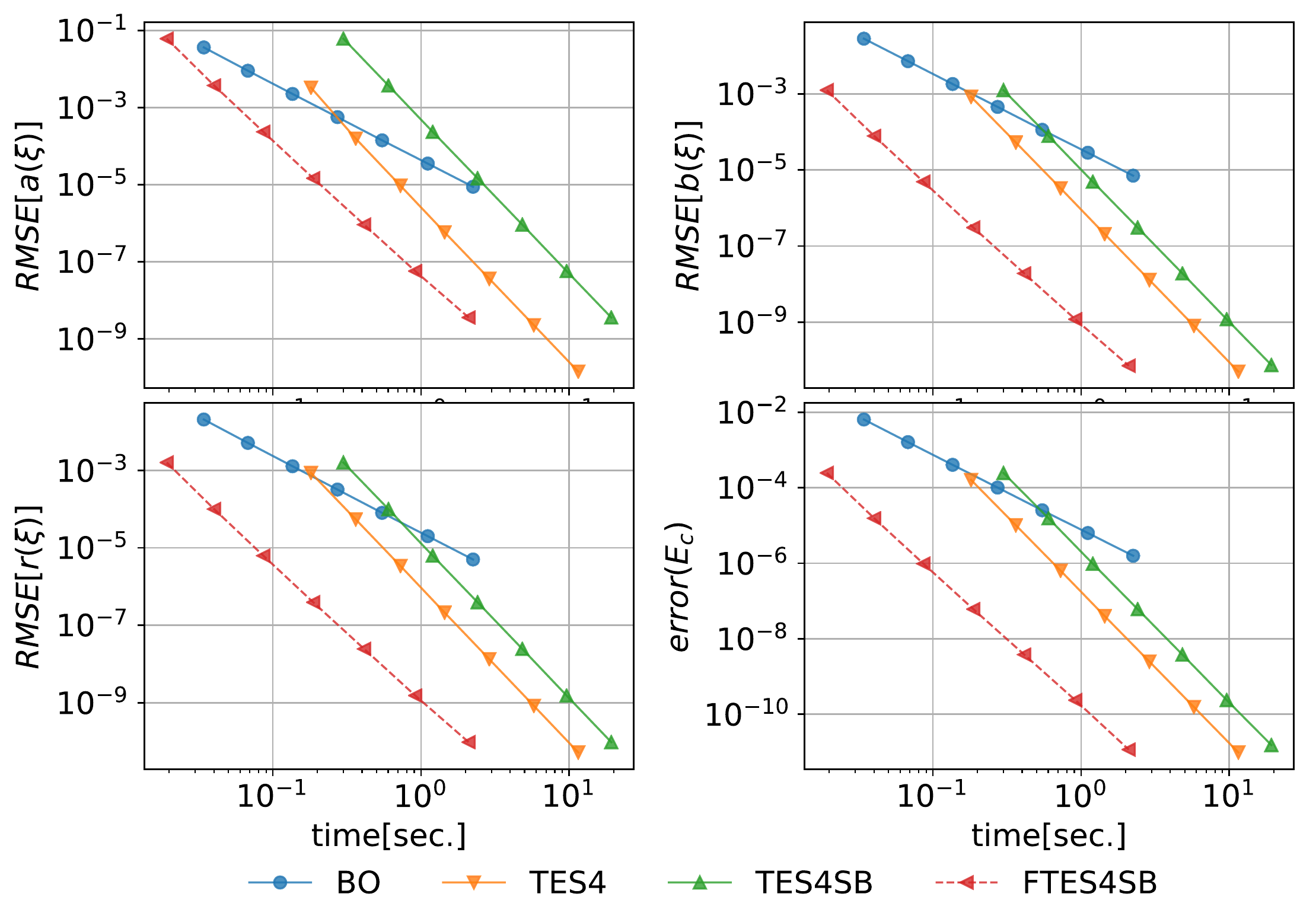}
\caption{The continuous spectrum errors depending on the execution time trade-off for the chirped hyperbolic secant~(\ref{Chirped}) in the case of normal dispersion~$\sigma=-1$.}
\label{fig:errorTimeSm1}
\end{figure}

\begin{figure}[htbp]
\centering
\includegraphics[width=0.6\linewidth]{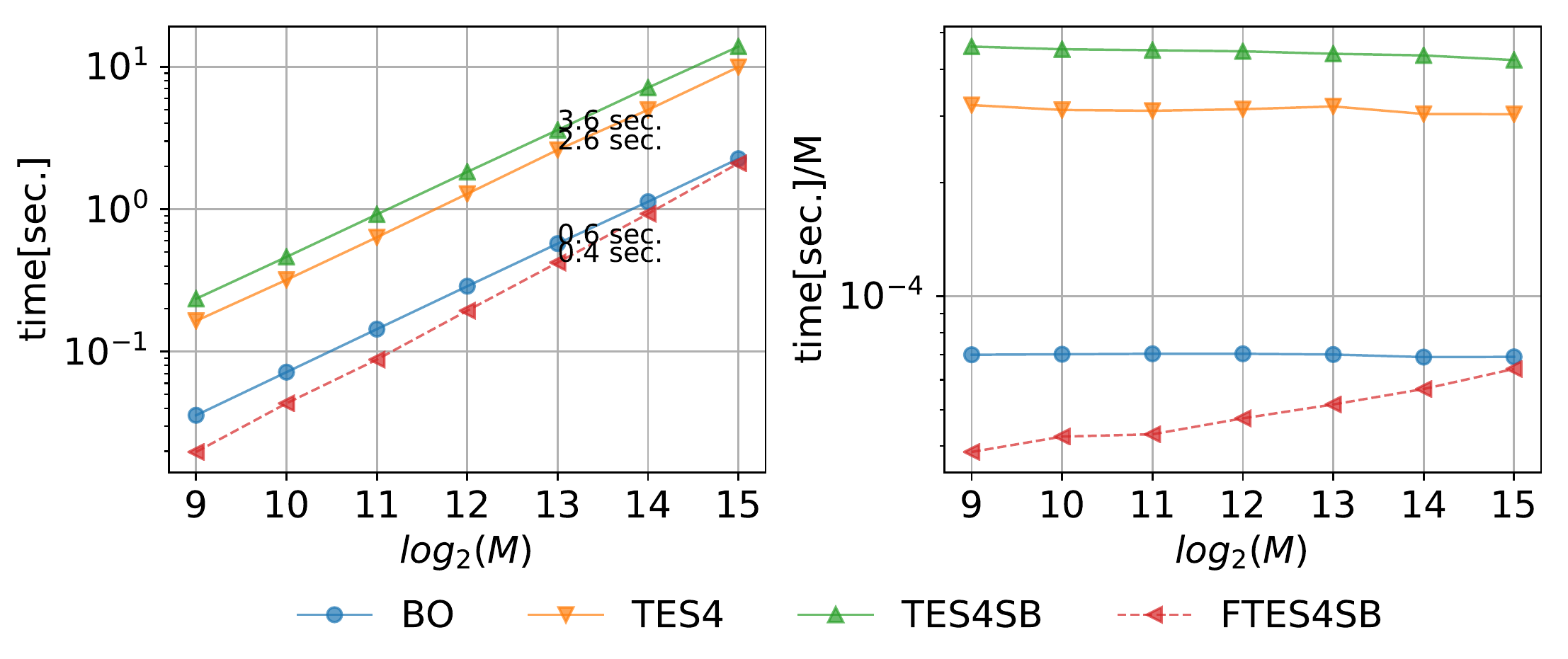}
\caption{Execution times for different algorithms in the case of anomalous dispersion~$\sigma=1$.}
\label{fig:timeS1}
\end{figure}

The efficiency of the schemes is compared in figures~\ref{fig:errorTimeS1} and \ref{fig:errorTimeSm1}. The fast variant of the proposed algorithm FTES4SB demonstrated the best speed when getting the desired error value across all considered schemes for both signs of dispersion. Of course, due to an asymptotic complexity of fast methods~\cite{Wahls2013}, one can determine the temporal grid size~$M$ for a fixed number of spectral parameter values~$N$ when the speed and efficiency of the fast scheme FTES4SB become comparable with conventional algorithms, which is demonstrated by Figure~\ref{fig:timeS1}. We should note here that the execution times of all methods don't depend on the signs of dispersion.

\begin{figure}[htbp]
\centering
\includegraphics[width=0.6\linewidth]{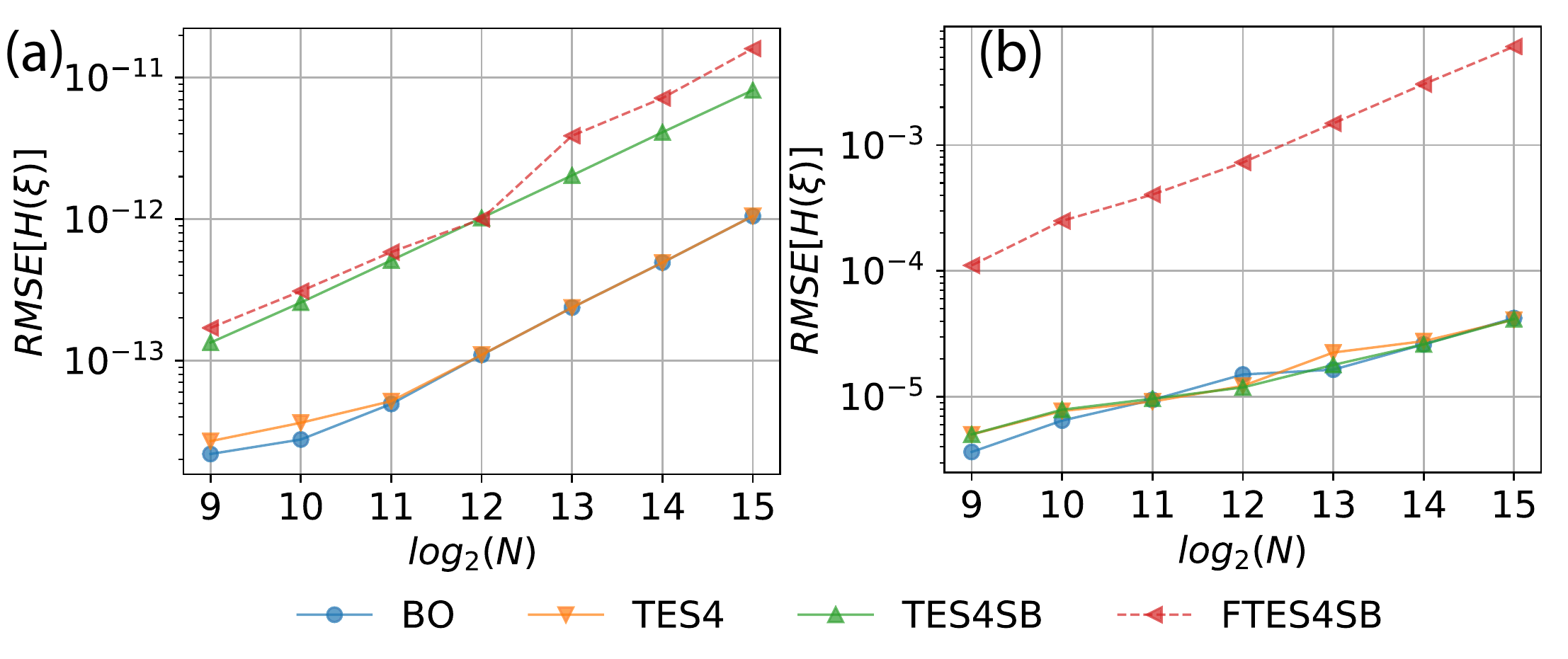}
\caption{Maximum value of the quadratic invariant~$H$ conservation error for anomalous dispersion~$\sigma=1$~(a) and normal dispersion~$\sigma=-1$~(b).}
\label{fig:invariantS1}
\end{figure}

\begin{figure}[htbp]
\centering
\includegraphics[width=0.7\linewidth]{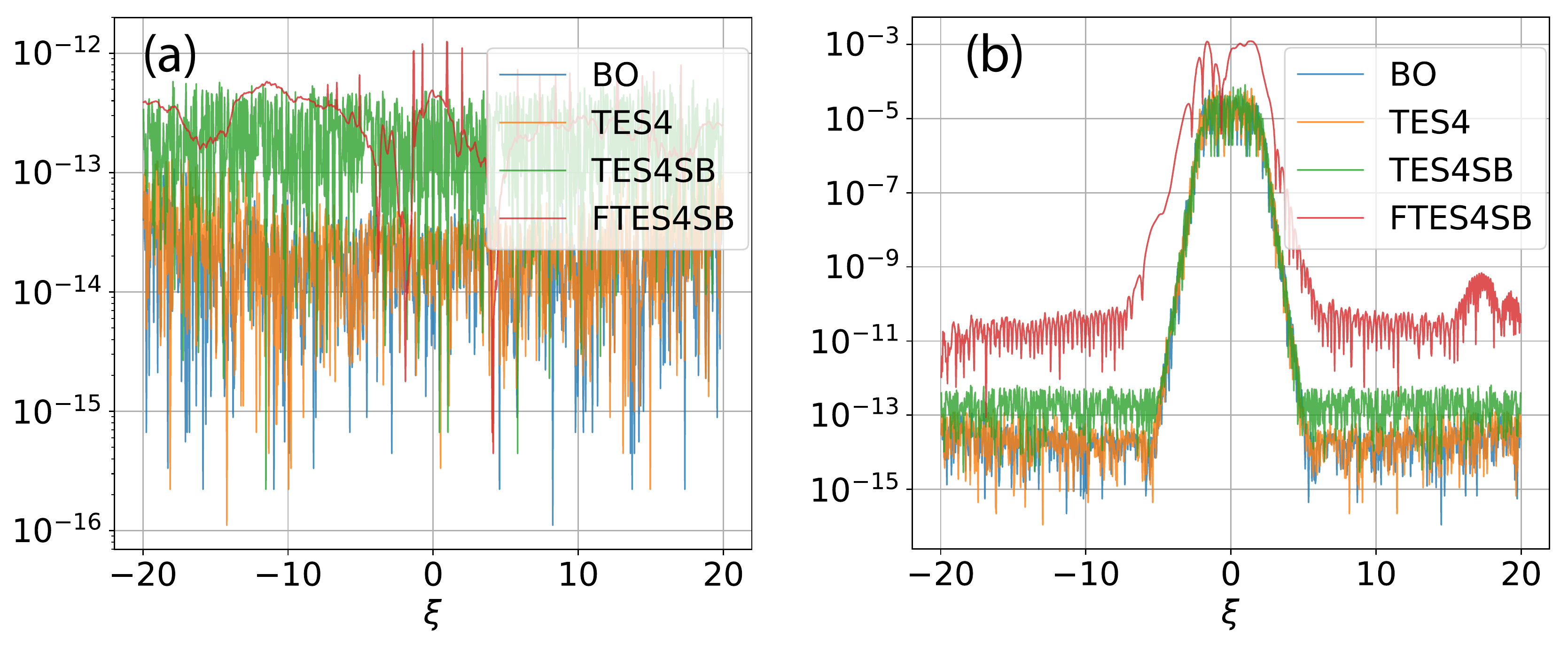}
\caption{Quadratic invariant conservation error~$|1 - |\psi_1|^2-\sigma|\psi_2|^2|$ depending on the spectral parameter~$\xi$: (a)~$\sigma=1$, (b)~$\sigma=-1$.}
\label{fig:invariant}
\end{figure}

The conservation properties of the schemes are considered in Figures~\ref{fig:invariantS1} and \ref{fig:invariant}. All algorithms demonstrated good conservation of the quadratic invariant~$H$ for the anomalous dispersion, but in case of normal dispersion, an error sufficiently increases. This is caused by the subtraction of large modulo quantities. All conventional schemes are comparable in the magnitude of the error. The accuracy of the proposed scheme in a fast variant (FTES4SB) reaches close value, though the fast computational technique caused an increase in error by two orders of magnitude for the normal dispersion.

\section{Conclusion}

In conclusion, we have developed a new multi-exponential scheme based on our three-exponential scheme and Suzuki decomposition, which allows fast computation and conserves the quadratic invariant for the real spectral parameter. The scheme consists of 13 matrix exponentials and has the 4th order of approximation. Also, it works for uniform grids, what together with the quadratic invariant conservation makes the proposed scheme attractive for telecommunication problems.

\section*{Funding}
Russian Science Foundation (RSF) (17-72-30006).

\bibliographystyle{unsrt}  
\bibliography{sample,references}  %%% Remove comment to use the external .bib file (using bibtex).
%%% and comment out the ``thebibliography'' section.

\end{document}